\documentclass[12pt,a4paper,twoside]{article}
\usepackage[utf8]{inputenc}
\usepackage[T1]{fontenc}
\usepackage{hyphsubst}
\usepackage[margin=2.5cm]{geometry}
\usepackage{lmodern}
\usepackage[square, comma, numbers, sort&compress, round]{natbib}
\usepackage[english]{babel}
\usepackage{bef_aras_nothm}
\usepackage{amsmath,amsthm,amssymb}
\usepackage{fixmath}
\usepackage{upgreek}
\usepackage{enumitem}
\usepackage{amsfonts}
\usepackage{microtype}
\usepackage{float}
\usepackage{mathrsfs}
\usepackage{hyperref}
\hypersetup{
     colorlinks   = true,
     citecolor    = blue
}
\usepackage[nameinlink]{cleveref} 
\crefname{figure}{\figurename}{\figurename}  
\usepackage{tikz}
\usepackage{verbatim}
\usepackage{afterpage}
\usepackage{glossaries}

\usepackage[headsepline]{scrlayer-scrpage}
\clearpairofpagestyles
\automark[subsection]{section}

\usetikzlibrary{arrows,calc}

\newglossary{symbols}{sym}{sbl}{Symbolverzeichnis}
\makeglossaries

\tikzset{
>=stealth',
help lines/.style={dashed, thick},
axis/.style={<->},
important line/.style={thick},
connection/.style={thick, dotted},
}

\theoremstyle{plain}
\newtheorem{thm}{Theorem}[section]
\newtheorem{lem}[thm]{Lemma}

\newtheorem{cor}[thm]{Corollary}
\theoremstyle{definition}

\newtheorem{rem}[thm]{Remark}

\newtheorem*{asA}{Assumption A.}
\newtheorem*{asB}{Assumption B.}

\newcommand{\Hi}{\mathscr{H}}
\newcommand{\F}{\mathfrak{F}}

\makeatletter
\g@addto@macro{\thm@space@setup}{\thm@headpunct{}}
\makeatother

\numberwithin{equation}{section}

\begin{document}

\newpage
\setcounter{page}{1}
\clearpairofpagestyles
\ohead{\pagemark}
\ihead{\headmark}

\title{Well-posedness of a fully nonlinear evolution inclusion of second order}

\author{Aras Bacho\footnotemark[2]}

\date{}
\maketitle

\footnotetext[1]{Ludwig-Maximilians-Universität München, Mathematisches Institut, Theresienstr. 39, 80333 München, Germany.}

\begin{abstract}
We consider the well-posedness of the abstract \textsc{Cauchy} problem for the doubly nonlinear evolution inclusion equation of second order given by
\begin{align*}
\begin{cases}
u''(t)+\partial \Psi(u'(t))+B(t,u(t))\ni f(t), &\quad t\in (0,T),\, T>0,\\
u(0)=u_0, \quad u'(0)=v_0 
\end{cases}
\end{align*}
where the function \( u \) takes values in a real separable \textsc{Hilbert} space, denoted by \(\mathscr{H}\). Here, \( u_0 \) lies in \(\mathscr{H}\), \( v_0 \) is in the intersection \(\overline{\mathrm{dom}(\partial \Psi)}\cap \mathrm{dom}(\Psi)\), and \( f \) belongs to \( \rmL^2(0,T;\mathscr{H}) \). The functional \(\Psi: \mathscr{H}\rightarrow (-\infty,+\infty] \) is assumed to be proper, lower semicontinuous, and convex. Additionally, the nonlinear operator \( B:[0,T]\times \mathscr{H}\rightarrow \mathscr{H} \) is assumed to satisfy either a global or a local \textsc{Lipschitz} condition. In the case where \( B \) satisfies a global Lipschitz condition, we can establish the existence and uniqueness of strong solutions \( u \) belonging to \( \rmH^2(0,T^*;\mathscr{H}) \). Furthermore, these solutions continuously depend on the data. We derive these results using the theory of nonlinear semigroups combined with the \textsc{Banach} fixed-point theorem. On the other hand, when \( B \) satisfies a local \textsc{Lipschitz} condition, we can guarantee the existence of strong local solutions.

\end{abstract}
\vspace*{1em} \textbf{Keywords} evolution inclusion of second order $ \cdot $ Well-posedness $ \cdot $ \textsc{Orlicz} space $ \cdot $ Fixed point argument $ \cdot $ Nonlinear semigroups $ \cdot $ Subdifferential operator 

\section{Introduction}
Throughout the article, let $(\Hi,\vert \cdot \vert,(\cdot,\cdot))$ be a real separable \textsc{Hilbert} space equipped with the inner product $(\cdot,\cdot)$ and the induced norm $\vert \cdot\vert :=(\cdot,\cdot)^{1/2}$ that we identify with its topological dual space $\Hi^*$ through the \textsc{Riesz} representation theorem, see, e.g., \textsc{Brézis} \cite[Theorem 4.11, p. 97]{Brez11FASS} . Then, we investigate the abstract \textsc{Cauchy} problem 
\begin{align} \label{eq:I}
\begin{cases}
u''(t)+\partial \Psi(u'(t))+B(t,u(t))\ni f(t) &\quad \text{for a.e. } t\in (0,T),\, T>0,\\
u(0)=u_0, \quad u'(0)=v_0 \tag{1}
\end{cases}
\end{align} where $\partial \Psi: \Hi\rightrightarrows \Hi$ is the \textsc{Fréchet} subdifferential of the functional $\Psi: \Hi\rightarrow (-\infty,+\infty]$, $B:[0,T]\times \mathscr{H}\rightarrow \mathscr{H}$ is a nonlinear operator, and $f:[0,T]\rightarrow \Hi$ is an external source.\\\\ We recall that for a proper, lower semicontinuous and convex functional $f:\Hi\rightarrow (-\infty,+\infty]$, the subdifferential of $f$ in $u\in \mathrm{dom}(f)$ is given by
\begin{align*}
\partial f(u)=\lbrace \xi\in \Hi : f(u)-f(v)\leq (\xi,u-v)\rbrace, 
\end{align*}  where $\mathrm{dom}(f)$ denotes the effective domain of $f$ defined by $\mathrm{dom}(f):= \lbrace u\in \Hi \mid f(u)<+\infty \rbrace$.
\subsection{Literature review}
Multivalued evolution equations (or evolution inclusions) of second order of the form
\begin{align}
 u''(t) + A(t)u'(t) + B(t)u(t) \ni f(t), \quad t \in (0,T),
\end{align}
have been studied in several cases and under various conditions by many authors. Primarily, two cases can be distinguished: In the first, the principal part of \( A \) is linear, and in the second, the principal part of the operator \( B \) is linear. 

Considering the first case, \textsc{B.} \cite{Bach22DNEI} examined the situation where \( A(t) = A_0 + A_1 \) is the sum of a linear, bounded, strongly positive, and symmetric operator \( A_0: V \to V^* \) defined on a reflexive and separable \textsc{Banach} space \( V \), and a multivalued perturbation \( A_1 = \partial \Psi_2 \) is given by the subdifferential of a proper, convex, and lower semicontinuous functional \( \Psi_2 \). The operator \( B(t) = \partial \calE_t + B(t, \cdot) \) is the sum of the subdifferential of a \( \lambda \)-convex functional \( \calE_t: V \to (-\infty, +\infty] \) with an effective domain within a reflexive and separable \textsc{Banach} space \( U \), and a strongly continuous perturbation \( B \). This equation is analyzed within a \textsc{Gelfand} triplet-type framework, without assuming \( U \) to be embedded in \( V \) or vice versa. Assuming that \( \calE_t \) satisfies a chain rule and \( \partial\calE_t \) adheres to a closedness condition, the existence of a strong solution satisfying an energy-dissipation inequality has been demonstrated using regularization techniques \cite{Bach23GMYR}. A more restricted result is presented in \cite{RosTho17CDIP}, where the rigorous assumptions prevent application to nonlinear elastodynamics. A related scenario for the single-valued case is explored in \cite{EmmSis13EESO}, where the operators satisfy an \textsc{Andrews--Ball} type condition, implying that \( (B + \lambda A): V \to V^* \) is monotone.
For the second case, where the principal component of \( B \) is linear, several authors have addressed the topic. In \cite{Bach21ONSA, BACHO2023126}, B. explored the \textsc{Cauchy} problem for the doubly nonlinear equation
\begin{align*}
u''(t) + \partial \Psi_{u(t)}(u'(t)) + \partial \calE_t(u(t)) + B(t,u(t),u'(t)) = f(t), \quad t \in (0,T),
\end{align*}
where \( \partial \Psi_{u(t)}(u'(t)) \) is multi-valued and nonlinear in both \( u' \) and \( u \), while the main component of \( \partial \calE_t(u(t)) \) is a linear, bounded, positive, and symmetric operator. The operator \( B \) may be seen as a strongly continuous perturbation of both \( \partial \Psi_{u(t)}(u'(t)) \) and \( \partial \calE_t(u(t)) \). For the scenario where \( A \) is linear, bounded, symmetric, and positive and \( B \) is maximally monotone, \textsc{Barbu} demonstrated the existence and uniqueness of strong solutions in \cite{Barb10NDMT}. The case where the acceleration term $u''$ is neglected, has been studied in \cite{BaEmMi19EREG, Bach21ONSA}.

In the context of single-valued operators, \textsc{Lions \& Strauss} in their seminal work \cite{LioStr65SNLE} established the well-posedness of the \textsc{Cauchy} problem for the doubly nonlinear evolution equation 
\begin{align*}
u''(t) + A(t,u(t),u'(t)) + B u(t) = f(t), \quad t \in (0,T),
\end{align*}
where \( B \) is a linear, self-adjoint, and unbounded operator. The operator \( A \) is nonlinear in \( u' \) but linear in \( u \), fulfilling a type of monotonicity condition.

\textsc{Emmrich \& Thalhammer} in \cite{EmSiTh15FDNS} proved the existence of solutions wherein, for every \( t \in [0,T] \), the operator \( A(t): V_A \to V_A^* \) is hemicontinuous and meets a particular growth condition so that \( A + \kappa I \) is both monotone and coercive. The operator \( B(t) = B_0 + C(t): V_B \to V_B^* \) represents the sum of a linear, bounded, symmetric, and strongly positive operator and a strongly continuous perturbation \( C(t) \), with consistent assumptions on \( V_A \) and \( V_B \).

In all aforementioned scenarios where the principal part of \( B \) is linear, operators \( A \) and \( B \) are typically defined on distinct spaces.

Doubly nonlinear evolution inclusions in which the leading parts of \(A\) and \(B\) are both nonlinear and which are defined on different spaces are not addressable within our framework. Nonetheless, for certain concrete problems, the existence of solutions has been demonstrated by leveraging the distinct structures of the operators. For instance, \textsc{Puhst} in \cite{Puhs15OEVF} established the existence of weak solutions assuming that the operators \(A\) and \(B\) are nonlocal. \textsc{Friedman \& Ne\v{c}as} in \cite{FriNec88SNWE} derived the existence of weak solutions based on the assumption that the operators are potential operators, which are twice differentiable with \textsc{Hessian} matrices that are uniformly positive definite and bounded. Meanwhile, both \textsc{Bul\'{\i}\v{c}ek, M\'{a}lek \& Rajagopal} \cite{BuMaRa12KVMG} and \textsc{Bul\'{\i}\v{c}ek, Kaplick\'{y} \& Steinhauer} \cite{BuKaSt1GKVM} confirmed the existence of weak solutions, presuming that the operators adhere to strong monotonicity, \textsc{Lipschitz}, and growth conditions. Under enhanced regularity conditions on the operators, these solutions have been verified as classical.

However, to the best of the author's knowledge, no abstract results exist for fully nonlinear evolution inclusions.

For further results on nonlinear evolution equations, we refer to \textsc{Leray} \cite{Lera53HDE},  \textsc{Dionne} \cite{Dion62SPCH},  \textsc{Emmrich \& Thalhammer} \cite{EmmTha10CTDD,EmmTha11DNEE}, \textsc{Emmrich, \& {\v{S}}i{\v{s}}ka} \cite{EmmSis11FDSO} including stochastic perturbations, \textsc{Emmrich,\textsc{{\v{S}}i{\v{s}}ka} \& Thalhammer} \cite{EmSiTh15FDNS} for a numerical analysis, \textsc{Emmrich, \textsc{{\v{S}}i{\v{s}}ka} \& Wr\'{o}blewska-Kami\'{n}ska} \cite{EmSiWr16ESOO} and \textsc{Ruf} \cite{Ruf17CFDS} for results on \textsc{Orlicz} spaces, and the monographs \textsc{Lions} \cite{Lion69QMRP}, \textsc{Lions \& Magenes} \cite[Chapitre 3.8]{LioMag72NHBV}, \textsc{Barbu} \cite[Chapter V]{Barb76NSDE}, \textsc{Wloka} \cite[Chapter V]{Wlok82PD}, \textsc{Zeidler} \cite[Chapter 33]{Zeid90NFA2b},  \textsc{Roub{\'\i}{\v{c}}ek} \cite[Chapter 11]{Roub13NPDE} and the references therein.\\
The list of literature presented in this section is not exhaustive.

\section{Assumptions and main results}
In this section, we establish the existence and uniqueness of solutions as well as the continuous dependence of the solution from the data specified below. Before we state the main result, we collect all the assumptions concerning the functional $\Psi$, the nonlinear operator $B$ as well as the external force $f$. 

\asA Let $\Psi: \Hi\rightarrow (-\infty,+\infty]$ be a  proper, lower semicontinuous, and convex functional such that $0\in \mathrm{dom}(\partial \Psi):=\lbrace v\in \Hi \mid \partial \Psi(u)\neq \emptyset \rbrace$.

\asB $ $
\begin{itemize}
\item[$i)$] Let $B:[0,T]\times \Hi \rightarrow \Hi$ satisfy the following local \textsc{Lipschitz} condition: For all $R>0$, there exists a function $\alpha_R \in \rmL^2(0,T;\mathbb{R}_0^+)$ such that
\begin{align}\label{ass1}
\vert B(t,u)-B(t,v)\vert \leq \alpha_R(t)\vert u-v\vert \quad\text{for all }u,v\in B_\Hi(0,R),
\end{align} and almost all $t\in[0,T]$, where $B_\Hi(0,R)$ denotes the closed ball in $\Hi$ with radius $R>0$ and center $0\in \Hi$.  Furthermore, there exists a function $g\in \rmL^2(0,T;\mathbb{R}_0^+)$, such that 
\begin{align*}
\vert B(t,0)\vert\leq g(t)\quad\text{for a.a. }t\in[0,T].
\end{align*}
\item[$ii)$] For all strongly measurable $v:[0,T]\rightarrow \Hi$, the map $t\mapsto B(t,v(t))$ is strongly measurable as a mapping from  $[0,T]$ to $\Hi$.
\end{itemize}

\begin{rem} From Assumption A, it follows that the subdifferential $\partial \Psi$ is a maximal monotone operator in the sense of \textsc{Brézis}, see \cite{Brez73OMMS}. In particular, for all $u,v\in \mathrm{dom}(\partial \Psi)$, there holds 
\begin{align} \label{eq:monotone}
    0\leq  (\eta -\xi, u-v) \quad \text{for all } \eta\in \partial \Psi(u),\, \xi \in \partial \Psi(v).
\end{align}

\end{rem}

\begin{rem} The condition $0\in \mathrm{dom}(\partial \Psi):=\lbrace v\in \Hi \mid \partial \Psi(u)\neq \emptyset \rbrace$ in Assumption A could in fact be replaced by the more general condition $\mathrm{dom}(\partial \Psi)\neq \emptyset$.
\end{rem}

\begin{rem} 
From the \textsc{Lipschitz} continuity, the square-integrability of $\alpha$, as well as Assumption B, we infer that the map $t\mapsto B(t,u)$ is in $\rmL^2(0,T;\Hi)$ for all $u\in \Hi$.
\end{rem} 

Having collected all assumptions, we are in the position to state the main result.

\begin{thm}\label{main} Let $\Psi: \Hi\rightarrow (-\infty,+\infty]$ and $B: [0,T]\times \Hi\rightarrow \Hi$ be given and satisfy Assumption A and Assumption B. Then, for every initial values $u_0\in \Hi, v_0\in \overline{\mathrm{dom}(\partial \Psi)}\cap \mathrm{dom}(\Psi)$ and every external source term $f\in \rmL^2(0,T;\Hi)$, there exists a unique local in time strong solution to \eqref{eq:I}, i.e., there exits $\tilde{T}>0$ and functions $u \in \rmH^2(0,\tilde{T};\Hi)$ and a $\eta:[0,T]\rightarrow \Hi$ strongly measurable such that
\begin{align}\label{eq:IIa}
u''(t)+\eta(t)+B(t,u(t))= f(t) \quad &\text{for a.e. } t\in (0,\tilde{T}),\\ \label{eq:IIb}
\eta(t)\in \partial \Psi(u'(t))  \quad &\text{for a.e. } t\in (0,\tilde{T}),
\end{align} and the initial conditions $u(0)=u_0$ and $u'(0)=v_0$ are fulfilled. If $B$ satisfies the global \textsc{Lipschitz} condition: There exists a function $\alpha\in \rmL^2(0,T)$ such that
\begin{align*}
\vert B(t,v)-B(t,w)\vert \leq \alpha(t)\vert v-w\vert  \quad \text{for all }v,w\in \Hi,
\end{align*} and almost every $t\in(0,T)$, then there exists a unique global solution $u\in \rmH^2(0,T;\Hi)$ to (1). Furthermore, the solution depends continuously on the data, i.e., let $u_1$ and $u_2$ be the solution to (1) associated with the data $(f,u_0^1,v_0^1)$ and $(g,u_0^2,v_0^2)$ from $\rmL^2(0,T;\Hi)\times \Hi\times \overline{\mathrm{dom}(\partial \Psi)}\cap \mathrm{dom}(\Psi)$, respectively. Then, there exists a constant $M>0$ such that 
\begin{align}\label{stab}
\Vert u_1-u_2\Vert_{\rmC([0,T];\Hi)}^2\leq M e^{\int_0^T\alpha(t) d t}\left( \vert u_0^1-u_0^2\vert^2+\vert v_0^1-v_0^2\vert^2+\Vert f-g\Vert^2_{\rmL^2(0,T;\Hi)}\right).
\end{align}
\end{thm}

\begin{proof}
 The main idea consists in rewriting the evolution equation as two coupled first order evolutions equations 
\begin{align} \label{eq:IIIa}
v'(t)+\partial \Psi(v(t))+B(t,u(t))\ni f(t) \quad &\text{for a.e. } t\in (0,T),\\ \label{eq:IIIb}
u'(t) =v(t)  \quad &\text{for all } t\in (0,T),\\
u(0)=u_0, \quad v(0)=v_0 \notag
\end{align}
and considering for fixed $u\in \rmC([0,T];\Hi)$ the auxiliary problem 
\begin{align} \label{eq:IV}
\begin{cases}
\tilde{v}'(t)+\partial \Psi(\tilde{v}(t))\ni\tilde{f}(t) \quad &\text{for a.e. } t\in (0,T),\\
\tilde{v}(0)=v_0,
\end{cases}
\end{align} where $\tilde{f}(t)= f(t)-B(t,u(t)),\, t\in[0,T]$. We notice that since $u$ is continuous and Assumption B holds, $\tilde{f}\in \rmL^2(0,T;\Hi)$ is ensured. \\ Then, based on the theory of nonlinear semigroups, existence and uniqueness of strong solutions $\tilde{v}\in \rmH^1(0,T;\Hi)$ for the auxiliary problem \eqref{eq:IV} such that the differential inclusion in \eqref{eq:IV} is satisfied almost everywhere in $(0,T)$ and $\Psi(u)\in \rmW^{1,1}(0,T)$ for every initial value $v_0\in \overline{\mathrm{dom}(\partial \Psi)}\cap \mathrm{dom}(\Psi)$ is well known, see for instance \textsc{Brézis} \cite[Proposition 3.12, p. 106]{Brez73OMMS} or \textsc{Barbu} \cite[Theorem 4.11, p. 158]{Barb10NDMT}. Furthermore, there exists a measurable selection $\xi\in \rmL^2(0,T;\Hi)$ such that $\xi(t)\in \partial\Psi(\tilde{v}(t))$ for almost every $t\in (0,T)$.\\
 Denoting with $J:\rmC([0,T];\Hi)\rightarrow \rmH^1([0,T];\Hi)$ the solution operator which maps the function $u\mapsto J(u)$ to the unique solution of \eqref{eq:IV},  we obtain for $u_1,u_2\in B_{ \rmC(0,T;\Hi)}(u_0,R)$ for fixed $R>0$, the inequality
\begin{align*}
\vert J(u_1)(t)-J(u_2)(t)\vert^2 & \leq 2\int_0^t \vert B(r,u_1(r))-B(r,u_2(r))\vert \vert J(u_1)(r)-J(u_2)(r)\vert dr\\
&\leq \int_0^t \vert J(u_1)(r)-J(u_2)(r)\vert^2 dr+ \int_0^t \alpha_R(r)\vert u_1(r)-u_2(r)\vert^2 dr
\end{align*} for all $0\leq t\leq T$.
Then, with \textsc{Gronwall}'s lemma, we obtain
\begin{align}\label{eq:V}
\vert J(u_1)(t)-J(u_2)(t)\vert^2  \leq e^t \int_0^t \alpha_R(r)\vert u_1(r)-u_2(r)\vert^2 dr \quad \text{for all } t\in [0,T].
\end{align}
Then, in order to show the existence solutions to the initial value problem (1), it suffices to show that the map $\F: B_{\rmC([0,\tilde{T}];\Hi)}(u_0,R)\rightarrow B_{\rmC([0,\tilde{T}];\Hi)}(u_0,R)$ with 
\begin{align}
\F(u)(t):=u_0+\int_0^t J(u)(s)\dd s, \quad t\in[0,\tilde{T}]
\end{align} possesses a fix point for a time-point $0<\tilde{T}\leq T$, where $B_{\rmC([0,\tilde{T}];\Hi)}(u_0,R)$ denotes the closed ball in $\rmC([0,\tilde{T}];\Hi)$ of radius $R>0$ and center $u_0$ which can be seen as constant function in $\rmC([0,\tilde{T}];\Hi)$. As soon as existence of a fixed point $u$ of $\F$ is shown, it follows
\begin{align}\label{eq:VI}
u(t)=\F(u)(t)=u_0+\int_0^t J(u)(s)\dd s=u_0+\int_0^t \tilde{v}(s)\dd s, \quad t\in[0,\tilde{T}],
\end{align} i.e., relation \eqref{eq:IIIb} holds and it follows $u\in \rmH^2(0,\tilde{T};\Hi)$. Since the operator $J$ maps the fixed point to the unique solution of the auxiliary problem \eqref{eq:IV}, the differential inclusion \eqref{eq:IIIa} holds as well. Finally, taking into account that the initial conditions are also satisfied, we deduce the existence of a strong solution to \eqref{eq:I}. We notice that since the resolvent operator $J$ maps continuous functions into (absolutley) continuous functions, the operator $\F$ itself maps continuous functions into continuous functions.  \\\\\textbf{Uniqueness:}\\\\ Before showing the existence of strong solutions, we establish uniqueness of solutions on the whole interval $[0,T]$. For this, we assume there are two solutions $u_1,u_2\in \rmH^2(0,T;\Hi)$ of \eqref{eq:I} to the same initial data. Then, by making use of \eqref{eq:V} and \eqref{eq:VI}, we obtain
\begin{align*}
\sup_{s\in[0,t]}\vert u_1(t)-u_2(t) \vert^2 &\leq \sup_{s\in[0,t]} \vert \int_0^s  J(u_1)(\tau)-J(u_2)(\tau) \dd \tau \vert^2\\
&\leq \left(\int_0^t \vert J(u_1)(s)-J(u_2)(s)\vert \dd s\right)^2\\
&\leq \sqrt{t} \int_0^t \vert J(u_1)(s)-J(u_2)(s)\vert^2 \dd s\\
&\leq \sqrt{t} \int_0^t e^s\int_0^s  2\alpha_R(\tau) \vert u_1(\tau)-u_2(\tau) \vert^2 \dd \tau \dd s\\
&\leq \sqrt{T}e^T \int_0^t \Vert u_1-u_2 \Vert^2_{\rmC([0,s];\Hi)} \int_0^s  2\alpha(\tau) \dd \tau \dd s,
\end{align*} where $R:=\sup_{t\in[0,T]}(\vert u_1(t)\vert +\vert u_2(t)\vert)$. Defining $a(t):= \Vert u_1-u_2 \Vert^2_{\rmC([0,t];\Hi)}$ and $\lambda(t):=\int_0^s  2\alpha_R(\tau) \dd \tau$, there holds $a,\lambda \in L^{\infty}(0,T)$ with $\lambda\geq 0$ a.e. in $(0,T)$ such that
\begin{align*}
a(t)\leq \int_0^t \lambda(s) a(s)\dd s \quad \text{for all }t\in [0,T].
\end{align*} \textsc{Gronwall}s Lemma yields immediately $a\equiv 0$ on $[0,T]$ so that $u_1=u_2$. \\\\ \textbf{Existence:}\\\\ In order to prove existence of local solutions, we make use of the \textsc{Banach} fixed-point theorem which provides the existence of a (unique) solution on a possibly small time-interval, i.e., we show existence of local solutions. Then, by iterating this procedure and making sure that the time interval do not minimize in each iteration step, global solution can be constructed.  Therefore, we need to check that the conditions of the \textsc{Banach} fixed-point theorem are fulfilled. Primarily, we show that for fixed $R>0$ the map $\F: B_{\rmC([0,\tilde{T}];\Hi)}(u_0,R)\rightarrow B_{\rmC([0,\tilde{T}];\Hi)}(u_0,R)$ is well defined for sufficiently small $\tilde{T}>0$, i.e., it maps the closed ball in $\rmC([0,\tilde{T}];\Hi)$ of radius $R>0$ and center $u_0$ into itself. In order to do so, we need the following a priori estimate: Let $u\in \rmC([0,\tilde{T}];\Hi)$, $v=J(u)$, and $\xi\in \rmL^2(0,T;\Hi)$ such that $\xi(t)\in \partial\Psi(v(t))$ for almost every $t\in (0,T)$. Furthermore, by Assumption A, there exists $\eta \in \Psi(0)$. Then, we obtain
\begin{align*}
\frac{1}{2}\frac{\dd}{\dd t}\vert J(u)(t) \vert^2&=\frac{1}{2}\frac{\dd}{\dd t}\vert v(t) \vert^2\\
&=(v'(t),v(t))\\
&\leq (v'(t),v(t))+(\xi(t)-\eta,v(t)-0)\\
&=(f(t)-B(t,u(t))-\eta,v(t))\\
&\leq (\vert f(t)\vert+ \vert B(t,u(t))\vert +\vert \eta \vert) \vert v(t)\vert\\
&\leq (\vert f(t)\vert+ \alpha_R(t)\vert u(t)\vert+ g(t)+\vert \eta\vert)\vert J(u)(t)\vert\ \quad \text{for a.e. }t\in  [0,T],
\end{align*} where we have tested the auxiliary problem \eqref{eq:IV} with its unique solution $v=J(u)$, the first inequality follows from the monotonicity of the subdifferential $\partial \Psi$, and the following equality follows from the fact that $v$ solves the auxiliary problem \eqref{eq:IV}. \\ This again yields by \textsc{Gronwall}s Lemma
\begin{align}
\sup_{s\in[0,t]}\vert J(u)(t)\vert &\leq \vert v_0 \vert +\int_0^t 2 (\vert f(s)\vert+ \alpha(s)\vert u(s)\vert+ g(s)+\vert \eta\vert)\dd s \notag\\
&= C+\int_0^t 2(\vert f(s)\vert+ \alpha(s)\vert u(s)\vert+ g(s))\dd s \quad \text{for all }t\in[0,T].
\end{align} where we used the fact that $J(u)(0)=\tilde{v}(0)=v_0$ and defined $C:= (\vert v_0\vert+T\vert\eta \vert)$. Note that since $\eta\in \partial\Psi(0)$ and $v_0\in \Hi$ are fixed, the value $C$ is constant. Employing (2.6), we obtain for $u\in B_{\rmC([0,\tilde{T}];\Hi)}(u_0,R)$
\begin{align*}
\Vert \F(u)-u_0 \Vert_{\rmC([0,\tilde{T}];\Hi)} &\leq \int_0^{\tilde{T}} \vert J(u)(t)\vert \dd t\\
&\leq \int_0^{\tilde{T}} C+2\int_0^t (\vert f(s)\vert+ \alpha(s)\vert u(s)\vert+ g(s))\dd s \dd t\\
&\leq \tilde{T}(C+2(\Vert f\Vert_{L^1(0,T)}+\Vert \alpha\Vert_{L^1(0,T)} \Vert u\Vert_{\rmC([0,T];\Hi)}+\Vert g\Vert_{L^1(0,T)}))\\
&\leq \tilde{T}(C+ 2(\Vert f\Vert_{L^1(0,T)}+\Vert \alpha\Vert_{L^1(0,T)}(R+\vert u_0\vert)+\Vert g\Vert_{L^1(0,T)}))\\
&\leq R
\end{align*} with ${\tilde{T}}\leq T_1:=R(\vert v_0 \vert+2(\Vert f\Vert_{L^1(0,T)}+\Vert \alpha\Vert_{L^1(0,T)}(R+\vert u_0\vert)+\Vert g\Vert_{L^1(0,T)}))^{-1}>0$.
Second, we show that for sufficiently small ${\tilde{T}}>0$, the map $\F$ is also a contraction. Let $u,v\in B_{\rmC([0,{\tilde{T}}];\Hi)}(u_0,R)$. Then, employing inequality \eqref{eq:V}, we obtain
\begin{align*}
\Vert \F(u)-\F(v) \Vert_{\rmC([0,{\tilde{T}}];\Hi)}&=\sup_{t\in[0,{\tilde{T}}]}\vert \int_0^t J(u)(s)-J(v)(s) \dd s\vert\\
&\leq  \int_0^{\tilde{T}} \vert J(u)(s)-J(v)(s)\vert \dd s\\
&\leq \int_0^{\tilde{T}} e^{t/2} \left( \int_0^t \alpha_R(r)\vert u_1(r)-u_2(r)\vert^2 dr\right)^{1/2} \dd s \\
&\leq \int_0^{\tilde{T}} \Vert u-v \Vert_{\rmC([0,t];\Hi)} e^{t/2}\left( \int_0^s  \alpha(\tau)   \dd \tau\right)^{1/2} \dd s\\
&\leq \int_0^{\tilde{T}} \Vert u-v \Vert_{\rmC([0,t];\Hi)} e^{T/2}\Vert \alpha\Vert_{L^1(0,T)}^{1/2}  \dd \tau \dd s\\
&\leq {\tilde{T}}(e^{T/2}\Vert \alpha\Vert_{L^1(0,T)}^{1/2})\Vert u-v\Vert_{\rmC([0,{\tilde{T}}];\Hi)}\\
&\leq L \Vert u-v\Vert_{\rmC([0,{\tilde{T}}];\Hi)},
\end{align*} where $L:={\tilde{T}}e^{T/2}\Vert \alpha\Vert_{L^1(0,T)}^{1/2}<1$ for ${\tilde{T}}< T_2:= (e^{T/2}\Vert \alpha\Vert_{L^1(0,T)}^{1/2})^{-1}>0$. Thus, by the \textsc{Banach} fixed-point theorem, there exists a unique solution $u\in \rmC([0,\tilde{T}],\Hi)$ to (1) on the time interval $[0,\tilde{T}]$ with $0<\tilde{T}<\min \lbrace T_1,T_2 \rbrace$. We assumed here without loss of generality that $\alpha \neq 0 $ in $\rmL^2(0,{\tilde{T}})$, otherwise $B$ would be constant almost everywhere and the assertion would be trivial.\\ Now, there are two possibilities to show global existence of solutions in the case where $B$ satisfies the global \textsc{Lipschitz} condition. The first possibility is to show the boundedness of the derivative of a solutions on the whole interval, such that blow ups of not only the solution itself but also of its derivative in finite time can not occur. This would lead to an interval of existence independent of the initial values. Then, applying successively the \textsc{Banach} fixed point theorem to the new initial value problem where the initial values are determined by the solution of the previous step, so that this procedure would cover the whole interval. Another possibility is to define the operator $\F$ on the whole space $\rmC([0,T];\Hi)$ equipped with a norm equivalent to the standard one and employ again the \textsc{Banach} fixed point theorem, where we need the equivalent norm to ensure the contractivity of $\F$. We tackle the problem with the latter option and define the operator $\F: \rmC(0,T;\Hi)\rightarrow \rmC(0,T;\Hi)$ as in (2.5), where we equip the space $C(0,T;\Hi)$ with the norm $\Vert v \Vert_{\chi}:=\sup_{t\in [0,T]}e^{-\tilde{L}t}\vert v(t)\vert$ with $\tilde{L}=2\Vert \alpha \Vert_{L^1(0,T)}$. Since $\F$ is obviously a self map, it remains to show that $\F$ is a contraction:
\begin{align*}
\Vert \F(u)-\F(v) \Vert_{\chi}&=\sup_{t\in[0,T]}e^{-\tilde{L} t}\vert \int_0^t J(u)(s)-J(v)(s) \dd s\vert\\
&\leq  \sup_{t\in[0,T]} e^{-\tilde{L}t}\int_0^t \vert J(u)(s)-J(v)(s)\vert \dd s\\
&\leq  \sup_{t\in[0,T]} 2e^{-\tilde{L}t}\int_0^t \int_0^s \alpha(\tau)\vert u(\tau)-v(\tau)\vert \dd \tau \dd s\\
&\leq  \sup_{t\in[0,T]} 2e^{-\tilde{L}t}\int_0^t  \sup_{\tau\in[0,s]} \vert u(\tau)-v(\tau)\vert\int_0^s \alpha(\tau)  \dd \tau \dd s\\
&\leq  \sup_{t\in[0,T]} 2e^{-\tilde{L}t} \Vert \alpha\Vert_{L^1(0,T)}\int_0^t  \sup_{\tau\in[0,s]} e^{\tilde{L}\tau}e^{-\tilde{L}\tau}\vert u(\tau)-v(\tau)\vert \dd s \\
&\leq  \sup_{t\in[0,T]} 2e^{-\tilde{L}t} \Vert \alpha\Vert_{L^1(0,T)}\int_0^t   e^{\tilde{L}s} \dd s \Vert u-v\Vert_{\chi}\\
&= \sup_{t\in[0,T]} 2e^{-\tilde{L}t} \Vert \alpha\Vert_{L^1(0,T)}\frac{e^{\tilde{L}t}-1}{\tilde{L}} \Vert u-v\Vert_{\chi}\\
&= \sup_{t\in[0,T]} (1-e^{-\tilde{L}t}) \Vert u-v\Vert_{\chi}\\
&= (1-e^{-\tilde{L}T}) \Vert u-v\Vert_{\chi}\\
\end{align*} Therefore, the map $\F$ is a contraction on $\rmC([0,T],\Hi)$ and by the \textsc{Banach} fixed point theorem there exists a unique fixed point $u\in \rmC([0,T],\Hi)$ which is a solution to (1). \\\\\textbf{Stability:}\\\\ Finally, we want to show the continuous dependence of the solution from the data. Let $u_1$ and $u_2$ be the solution to (1) associated with $(f,u_0^1,v_0^1),(g,u_0^2,v_0^2)\in \rmL^2(0,T;\Hi)\times \Hi\times \overline{\mathrm{dom}(\Psi)}$, respectively. With the same reasoning as for \eqref{eq:V}, we derive with \textsc{Gronwall}'s lemma
\begin{align*}
\vert J(u_1)(t)-J(u_2)(t)\vert^2  \leq e^t\left(\vert v_0^1-v_0^2\vert^2+\Vert f-g\Vert^2_{\rmL^2(0,T;\Hi)}+ \int_0^t \alpha(r)\vert u_1(r)-u_2(r)\vert^2 dr\right) 
\end{align*} for all $t\in [0,T]$. Then, continuing as in the uniqueness part, we obtain
\begin{align*}
\sup_{s\in[0,t]}\vert u_1(s)-u_2(s) \vert^2 &\leq\vert u_0^1-u_0^1 \vert^2+\sup_{s\in[0,t]} \vert \int_0^s  J(u_1)(\tau)-J(u_2)(\tau) \dd \tau \vert^2\\
&\leq \vert u_0^1-u_0^1 \vert^2+ \left(\int_0^t \vert J(u_1)(s)-J(u_2)(s)\vert \dd s\right)^2\\
&\leq \vert u_0^1-u_0^1 \vert^2+ \sqrt{t} \int_0^t \vert J(u_1)(s)-J(u_2)(s)\vert^2 \dd s\\
&\leq \vert u_0^1-u_0^1 \vert^2+\sqrt{t} \int_0^t e^s\left(\vert v_0^1-v_0^2\vert^2+\Vert f-g\Vert^2_{\rmL^2(0,T;\Hi)}\right)\dd s\\
&\quad+\int_0^t e^s\int_0^s  2\alpha_R(\tau) \vert u_1(\tau)-u_2(\tau) \vert^2 \dd \tau \dd s\\
&\leq M\left(\vert u_0^1-u_0^1 \vert^2+\vert v_0^1-v_0^2\vert^2+\Vert f-g\Vert^2_{\rmL^2(0,T;\Hi)}\right)\\
&\quad+\sqrt{T}e^T \int_0^t \Vert u_1-u_2 \Vert^2_{\rmC([0,s];\Hi)} \int_0^s  2\alpha(\tau) \dd \tau \dd s,
\end{align*} for a constant $M>0$ independent of the data. Defining again $a(t):= \Vert u_1-u_2 \Vert^2_{\rmC([0,t];\Hi)}$ and $\lambda(t):=\int_0^s  2\alpha_R(\tau) \dd \tau$ as well as $b=M\left(\vert u_0^1-u_0^1 \vert^2+\vert v_0^1-v_0^2\vert^2+\Vert f-g\Vert^2_{\rmL^2(0,T;\Hi)}\right)$, there holds $a,\lambda \in L^{\infty}(0,T)$ with $\lambda\geq 0$ a.e. in $(0,T)$ such that
\begin{align*}
a(t)\leq b+\int_0^t \lambda(s) a(s)\dd s \quad \text{for all }t\in [0,T].
\end{align*} Again, with the \textsc{Gronwall} lemma, we obtain the desired estimate.

\end{proof} 
\begin{cor} In the case, when $B$ satisfies the local \textsc{Lipschitz} condition, there exists a maximal solution, i.e., there exists a time interval $[0,\tilde{T})\subset [0,T]$ and a function $u$ such that for each compact subinterval $[0,S] \subset [0,\tilde{T}]$ there holds $u\in \rmH^2(0,S,\Hi)$ and $u$ solves problem (1) pointwise almost every on $(0,\tilde{T})$. Furthermore, for every sequence $(t_n)\subset [0,\tilde{T})$ with $t_n \nearrow \tilde{T}$ as $n\rightarrow \infty$, there holds $\vert u(t_n)\vert \rightarrow +\infty$ as $n\rightarrow \infty$.
\end{cor}

\begin{rem} We notice that we did not impose any compactness assumption neither on the sublevels of the dissipation potential $\Psi$ nor that the operator $B$ is a strongly continuous perturbation in order to show existence of solutions.
\end{rem}

\section{Example}
In this section, we want to apply the abstract result to concrete examples. Let, $\Omega\subset \mathbb{R}^d$ be a \textsc{Lipschitz} domain and $T>0$. We first consider the following initial-boundary value problem

\begin{align*}
\text{(P1)}
\begin{cases}
\partial_{tt} u(\xx,t)-\nabla_{\xx} \cdot \mathbf{p}(\xx,t) -\nabla_{\xx}\cdot H(t,\xx,\nabla_{\xx} u) = f(\xx,t) \quad &\text{in } \Omega\times (0,T),\\
\mathbf{p}(\xx,t)\in \partial_{\zz} \psi(\xx, \nabla_{\xx} \partial_t u(\xx,t)) \quad &\text{a.e. in }  \Omega\times (0,T),\\
u(\xx,0)\,\,=u_0(\xx) \quad &\text{on } \Omega,\\
u'(\xx,0)\,=v_0(\xx) \quad \,&\text{on } \Omega, \\
u(\xx,t)\,\,\,=0 \quad  \qquad &\text{on } \partial \Omega\times[0,T],\\
\frac{\partial u}{\partial \nu}(\xx,t)=0 \quad  \qquad &\text{on } \partial \Omega\times[0,T]
\end{cases}
\end{align*} where $f:\Omega \times [0,T]\rightarrow \mathbb{R}, b:\Omega \times [0,T]\times \mathbb{R}\rightarrow \mathbb{R}$ and $\psi:\Omega\times \mathbb{R}^d\rightarrow (-\infty,+\infty]$ are measurable functions satisfying the following conditions:
\begin{enumerate} [label=(\thesection.\alph*), leftmargin=3.2em]  \label{eq:psi.applic.II}
\item \label{2.a} The function $\psi:\Omega\times \mathbb{R}^d \rightarrow [0,+\infty]$ is a non-negative  \textsc{Carath\'{e}odory} function such that $\psi(\xx, \cdot)$ is a proper, lower semicontinuous, and convex, and $\psi(\xx,0)=0$ for almost every $\xx\in \Omega$.
\item \label{2.b} There exists a strictly increasing, convex, and lower semicontinuous function $\Phi:\mathbb{R}\rightarrow [0,+\infty]$ with 
\begin{align*}
\lim_{t \rightarrow +\infty}\frac{\Phi(t)}{t}=+\infty
\end{align*}
such that 
\begin{align*}
\Phi\left(\vert \xii \vert\right)\leq \psi(\xx,\xii)
\end{align*} $\text{for a.e. }\xx\in \Omega \text{ and all }\xii \in \mathbb{R}^d$. 
\item \label{2.c} The function $H:[0,T]\times \Omega\times \mathbb{R}^d \rightarrow \mathbb{R}^d$ is a \textsc{Carath\'{e}odory} function satisfying the following \textsc{Lipschitz} condition: there exist a constant $C_H>0$ such that for almost all $\xx \in \Omega$ and $t\in (0,T)$
\begin{align*}
\vert H(t,\xx,\xii_1)-H(t,\xx,\xii_2)\vert\leq C_H\vert \xii_1-\xii_2\vert \quad \text{ for all }\xii_1,\xii_2 \in \mathbb{R}^d.
\end{align*}
\item \label{2.d} There holds $f\in \rmL^2 (0,T;\rmH_0^1(\Omega))$.
\end{enumerate}
Then, we naturally choose $\Hi=\rmH_0^1(\Omega)$ and identify $\Hi$ with its dual space $\Hi^*=\rmH^{-1}(\Omega)$ by the \textsc{Riesz} isomporphism. Then, the functional $\Psi:\Hi\rightarrow \mathbb{R}$ and the operator $B:[0,T]\times \Hi \rightarrow \Hi$ are given by 
\begin{align*}
\Psi(v)=\int_{\Omega} \psi(\xx,\nabla v(\xx)) \dd \xx 
\end{align*} and
\begin{align*}
(B(t,u),w )= \int_{\Omega} H(t,\xx,\nabla u(\xx)) \cdot \nabla w(\xx) \dd \xx \quad \text{for all }w\in \Hi.
\end{align*} An important class examples for the function $\psi$ is given by the class of convex superlinear and anisotropic \textsc{Orlicz} functions, e.g., 
\begin{itemize}
\item[1.] $\psi(\xii)=\vert \xii\vert \log(1+\vert \xii \vert)$ with $\partial \psi(\xii)=\mathrm{Sgn}(\xii)\log(1+\vert \xii \vert)+\frac{\xii}{1+\vert \xii \vert}$,
\item[2.] $\psi(\xii)=\vert \xii \vert \exp(\vert \xii \vert)$ with $\partial \psi(\xii)=\mathrm{Sgn}(\xii)\exp(\vert \xii \vert)+\xii\exp(\vert \xii \vert)$,
\item[3.] $\psi(\xii)= \frac{1}{p}\vert \xii \vert^{p}+\vert \xii \vert$ for $p> 1$ with $\partial \psi(\xii)=\mathrm{Sgn}(\xii)+\xii\vert \xii \vert^{p-2}$,
\item[4.] $\psi(\xii)= \exp(\frac{1}{p}\vert \xii \vert^p+\vert \xii \vert)$ for $p>1$ with $\partial \psi(\xii)=\exp(\frac{1}{p}\vert \xii \vert^p+\vert \xii \vert)(\mathrm{Sgn}(\xii)+\xii\vert \xii \vert^{p-2})$,
\item[5.] $\psi(\xii)= \exp(\vert \xii\vert \log(1+\vert \xii \vert))$ with $\partial \psi(\xii)=\exp(\vert \xii\vert \log(1+\vert \xii \vert))(\mathrm{Sgn}(\xii)\log(1+\vert \xii \vert)+\frac{\xii}{1+\vert \xii \vert})$,
\end{itemize} and so on, where $\mathrm{Sgn}:\mathbb{R}^d\rightarrow \mathbb{R}^d$ denotes the multi-valued and multi-dimensional sign function defined by
\begin{align*}
\mathrm{Sgn}(\xii)=
\begin{cases}
B(0,1), &\text{if }\xii=0,\\
\frac{\xii}{\vert \xii\vert},  &\text{ otherwise}.
\end{cases}
\end{align*} For more example, see also \cite{EmmWro13CQPE,CaZh09EUWP,Dona74IONP,GwSw11PEAO,Zeid85NFA3}. In that case, the effective domain of the functional $\Psi$ is given by an anisotropic \textsc{Orlicz} space which is, in general, neither  reflexive nor separable which makes the analysis in general difficult. 
For the function $H$, we can have $H(\xii)=\xii (\vert \xii\vert^2+1)^{\frac{q-1}{2}} $ and $H(\xii)=\xii (\vert \xii\vert+1)^{q-1} $ with $1\leq q\leq 2$, $H(\xii)=\mathrm{tanh}(\vert \xii\vert^2)$ or any other globally \textsc{Lipschitz} function $H$.\\

In the following lemma, we show that conditions of Theorem \ref{main}
are satisfied.
\begin{lem} Let the conditions \ref{2.a}-\ref{2.c} be fulfilled and assume that there exists $\tilde{v}\in \rmW^{1,\infty}(\Omega)$ such that $\Psi(\tilde{v})<+\infty$. Then, the functional $\Psi:\Hi\rightarrow [0.+\infty]$ is proper, lower semicontinuous, and convex and the operator $B;[0,T]\times \Hi \rightarrow \Hi$ is \textsc{Lipschitz} continuous. Furthermore, let $v\in \mathrm{dom}(\partial \Psi)$. Then, there holds $\zeta \in \partial \Psi(v)$ if and only if $\zeta=\nabla \cdot \mathbf{p}$ in the weak sense and $\mathbf{p}(\xx)\in \partial_{\zz}\psi(\xx,\nabla v(\xx))$ for almost every $\xx\in \Omega$.

\end{lem}
\begin{proof} It is easily checked that $\Psi$ is proper and convex on $\Hi$. We show that $\Psi$ is lower semicontinuous on $\Hi$. For that, we equivalently show that the sublevel sets $J_\alpha:=\lbrace v\in \Hi : \Psi(v)\leq \alpha\rbrace$ are closed for all $\alpha\in \mathbb{R}$. Let $\alpha\in \mathbb{R}$ and let $(v_n)_{n\in \mathbb{N}}\subset J_\alpha$ such that $v_n\rightarrow v$ in $\Hi$ as $n\rightarrow \infty$. By Condition \ref{2.b} and the \textsc{de la Vallée Poussin} criterion, we know that $(\nabla v_n)_{n\in \mathbb{N}}$ is weakly compact in $\rmL^1(\Omega)$ and therefore (up to a subsequence) $\nabla v_n\rightharpoonup \nabla v$ in $\rmL^1(\Omega)$. By uniqueness of the limit, we infer the convergence of the whole sequence. Now, choose a subsequence $n_k$ such that 
\begin{align*}
\lim_{k\rightarrow \infty} \int_\Omega \psi(\xx,\nabla v_{n_k}(\xx))\dd \xx =\liminf_{n\rightarrow \infty}\int_\Omega \psi(\xx,\nabla v_{n}(\xx))\dd \xx.
\end{align*} Then, by \textsc{Ekeland} and \textsc{Temam} \cite[Theorem 2.1, p. 243]{EkeTem76CAVP}, there holds 
\begin{align*}
\int_\Omega \psi(\xx,\nabla v(\xx))\dd \xx &\leq \liminf_{k\rightarrow \infty}\int_\Omega \psi(\xx,\nabla v_{n_k}(\xx))\dd \xx\\
& \leq \liminf_{n\rightarrow \infty}\int_\Omega \psi(\xx,\nabla v_{n}(\xx))\dd \xx,
\end{align*} which shows the lower semicontinuity of $\Psi$ on $\Hi$. \\ Now, we show the \textsc{Lipschitz} continuity of $B$:
\begin{align*}
(B(t,u)-B(t,v),w )&= \int_{\Omega}\left( b(t,\xx,u(\xx))- b(t,\xx,v(\xx)) \right) w(\xx) \dd \xx\\
&\leq \left(\int_{\Omega}\left\vert b(t,\xx,u(\xx))- b(t,\xx,v(\xx)) \right \vert^2 \dd \xx\right)^{1/2}\Vert w\Vert_{\rmL^2(\Omega)}\\
&\leq \left(\int_{\Omega} C^2_b\vert u(\xx) - v(\xx) \vert^2 \dd \xx\right)^{1/2}\Vert w\Vert_{\rmL^2(\Omega)}  \quad \text{for all }w\in \Hi
\end{align*} and hence the \textsc{Lipschitz} continuity of $B$. The last assertion follows from Proposition 5.1, p. 21, Proposition 5.7, p. 27, and Proposition 2.1, p. 271, in \textsc{Ekeland} and \textsc{Temam} \cite{EkeTem76CAVP} which completes the proof.
\end{proof}
We are now in the position, to state the main result.
\begin{thm}
Let the Conditions \ref{2.a}-\ref{2.d} be satisfied and assume that there exists $\tilde{v}\in \rmW^{1,\infty}(\Omega)$ such that $\Psi(\tilde{v})<+\infty$. Then, for every $u_0\in \rmL^2(\Omega) $ and  $v_0\in \overline{\mathrm{dom}(\partial \Psi)}\cap \mathrm{dom}(\Psi)$, there exists a unique weak solution $u\in \rmH^2(0,T,\Hi)$ to \textnormal{(P1)}. Furthermore, the stability estimate \ref{stab} is fulfilled.
\end{thm}
\begin{proof}
Since the assumptions are all fulfilled, this follows immediately from  Theorem \ref{main}.
\end{proof}
\begin{rem}
We notice again, that we did not impose any compactness condition on $\Psi$ but rely on the maximal monotonicity of its subdifferential. This implies that we do not handle the nonlinearity $B$ by compactness methods but a fixed point argument. Compare with \cite{EmSiWr16ESOT}, where a time discretization method has been employed to show with compactness arguments the existence of weak solutions in an \textsc{Orlicz} space.
\end{rem}

\begin{rem} Instead of the global \textsc{Lipschitz} condition \ref{2.c}, we can impose a local \textsc{Lipschitz} condition on $B$. Consequently, we obtain the local existence result given in Theorem \ref{main}. 
\end{rem}

Next, we study the following initial-boundary value problem

\begin{align*}
	\text{(P2)}
	\begin{cases}
		\partial_{tt} u(\xx,t)-\nabla_{\xx} \cdot \mathbf{p}(\xx,t) +b(t,\xx,\nabla_{\xx} u) = f(\xx,t) \quad &\text{in } \Omega\times (0,T),\\
		\mathbf{p}(\xx,t)\in \partial_{\zz} \psi(\xx, \nabla_{\xx} \partial_t u(\xx,t)) \quad &\text{a.e. in }  \Omega\times (0,T),\\
		u(\xx,0)\,\,=u_0(\xx) \quad &\text{on } \Omega,\\
		u'(\xx,0)\,=v_0(\xx) \quad \,&\text{on } \Omega, \\
		u(\xx,t)\,\,\,=0 \quad  \qquad &\text{on } \partial \Omega\times[0,T],\\
		\frac{\partial u}{\partial \nu}(\xx,t)=0 \quad  \qquad &\text{on } \partial \Omega\times[0,T]
	\end{cases}
\end{align*} where $f:\Omega \times [0,T]\rightarrow \mathbb{R}, b:\Omega \times [0,T]\times \mathbb{R}\rightarrow \mathbb{R}$ and $\psi:\Omega\times \mathbb{R}^d\rightarrow (-\infty,+\infty]$ are measurable functions satisfying the following conditions:
\begin{enumerate} [label=(\thesection.\alph*), leftmargin=3.2em]  \label{eq:psi.applic.II}
	\item \label{2.a} The function $\psi:\Omega\times \mathbb{R}^d \rightarrow [0,+\infty]$ satisfies Conditions (3.a) and (3.b) as in the first example. 
	\item \label{2.b} There exists a strictly increasing, convex, and lower semicontinuous function $\Phi:\mathbb{R}\rightarrow [0,+\infty]$ with 
	\begin{align*}
		\lim_{t \rightarrow +\infty}\frac{\Phi(t)}{t}=+\infty
	\end{align*}
	such that 
	\begin{align*}
		\Phi\left(\vert \xii \vert\right)\leq \psi(\xx,\xii)
	\end{align*} $\text{for a.e. }\xx\in \Omega \text{ and all }\xii \in \mathbb{R}^d$. 
	\item \label{2.c} The function $b:[0,T]\times \Omega\rightarrow \mathbb{R}$ is a \textsc{Carath\'{e}odory} function satisfying the following \textsc{Lipschitz} condition: there exist a constant $C_b>0$ such that
	\begin{align*}
		\vert b(t,\xx,u)-b(t,\xx,v)\vert\leq C_b\vert u-v\vert \quad \text{ for all }u,v\in \mathbb{R}.
	\end{align*}
	\item \label{2.d} There holds $f\in \rmL^2(\Omega\times (0,T))$.
\end{enumerate}
Then, we naturally choose $\Hi=\rmL^2(\Omega)$. Then, the functional $\Psi:\Hi\rightarrow \mathbb{R}$ and operator $B:[0,T]\times \Hi \rightarrow \Hi$ are given by 
\begin{align*}
	\Psi(v)=\int_{\Omega} \psi(\xx,\nabla v(\xx)) \dd \xx 
\end{align*} and
\begin{align*}
	(B(t,u),w )= \int_{\Omega} b(t,\xx,u(\xx)) w(\xx) \dd \xx \quad \text{for all }w\in \Hi.
\end{align*} An important class examples for the function $\psi$ is given by the class of convex superlinear and anisotropic \textsc{Orlicz} functions, e.g., 
\begin{itemize}
	\item[1] $\psi(\xii)=\vert \xii\vert \log(1+\vert \xii \vert)$ with $\partial \psi(\xii)=\mathrm{Sgn}(\xii)\log(1+\vert \xii \vert)+\frac{\xii}{1+\vert \xii \vert}$,
	\item[2] $\psi(\xii)=\vert \xii \vert \exp(\vert \xii \vert)$ with $\partial \psi(\xii)=\mathrm{Sgn}(\xii)\exp(\vert \xii \vert)+\xii\exp(\vert \xii \vert)$,
	\item[3] $\psi(\xii)= \frac{1}{p}\vert \xii \vert^{p}+\vert \xii \vert$ for $p> 1$ with $\partial \psi(\xii)=\mathrm{Sgn}(\xii)+\xii\vert \xii \vert^{p-2}$,
	\item[4] $\psi(\xii)= \exp(\frac{1}{p}\vert \xii \vert^p+\vert \xii \vert)$ for $p>1$ with $\partial \psi(\xii)=\exp(\frac{1}{p}\vert \xii \vert^p+\vert \xii \vert)(\mathrm{Sgn}(\xii)+\xii\vert \xii \vert^{p-2})$,
	\item[5] $\psi(\xii)= \exp(\vert \xii\vert \log(1+\vert \xii \vert))$ with $\partial \psi(\xii)=\exp(\vert \xii\vert \log(1+\vert \xii \vert))(\mathrm{Sgn}(\xii)\log(1+\vert \xii \vert)+\frac{\xii}{1+\vert \xii \vert})$,
\end{itemize} and so on, where $\mathrm{Sgn}:\mathbb{R}^d\rightarrow \mathbb{R}^d$ denotes the multi-valued and multi-dimensional sign function defined by
\begin{align*}
	\mathrm{Sgn}(\xii)=
	\begin{cases}
		B(0,1), &\text{if }\xii=0,\\
		\frac{\xii}{\vert \xii\vert},  &\text{ otherwise}.
	\end{cases}
\end{align*} For more example, see also \cite{EmmWro13CQPE,CaZh09EUWP,Dona74IONP,GwSw11PEAO,Zeid85NFA3}. In that case, the effective domain of the functional $\Psi$ is given by an anisotropic \textsc{Orlicz} space which is, in general, neither  reflexive nor separable which makes the analysis in general difficult.
In the following lemma, we show that conditions of Theorem \ref{main}
are satisfied.
\begin{lem} Let the conditions \ref{2.a}-\ref{2.c} be fulfilled and assume that there exists $\tilde{v}\in \rmW^{1,\infty}(\Omega)$ such that $\Psi(\tilde{v})<+\infty$. Then, the functional $\Psi:\Hi\rightarrow [0.+\infty]$ is proper, lower semicontinuous, and convex and the operator $B;[0,T]\times \Hi \rightarrow \Hi$ is \textsc{Lipschitz} continuous. Furthermore, let $v\in \mathrm{dom}(\partial \Psi)$. Then, there holds $\zeta \in \partial \Psi(v)$ if and only if $\zeta=\nabla \cdot \mathbf{p}$ and $\mathbf{p}(\xx)\in \partial_{\zz}\psi(\xx,\nabla v(\xx))$ for almost every $\xx\in \Omega$.

\end{lem}
\begin{proof} It is easily checked that $\Psi$ is proper and convex on $\Hi$. We show that $\Psi$ is lower semicontinuous on $\Hi$. For that, we equivalently show that the sublevel sets $J_\alpha:=\lbrace v\in \Hi : \Psi(v)\leq \alpha\rbrace$ are closed for all $\alpha\in \mathbb{R}$. Let $\alpha\in \mathbb{R}$ and let $(v_n)_{n\in \mathbb{N}}\subset J_\alpha$ such that $v_n\rightarrow v$ in $\Hi$ as $n\rightarrow \infty$. By Condition \ref{2.b} and the \textsc{de la Vallée Poussin} criterion, we know that $(\nabla v_n)_{n\in \mathbb{N}}$ is weakly compact in $\rmL^1(\Omega)$ and therefore (up to a subsequence) $\nabla v_n\rightharpoonup \nabla v$ in $\rmL^1(\Omega)$. By uniqueness of the limit, we infer the convergence of the whole sequence. Now, choose a subsequence $n_k$ such that 
	\begin{align*}
		\lim_{k\rightarrow \infty} \int_\Omega \psi(\xx,\nabla v_{n_k}(\xx))\dd \xx =\liminf_{n\rightarrow \infty}\int_\Omega \psi(\xx,\nabla v_{n}(\xx))\dd \xx.
	\end{align*} Then, by \textsc{Ekeland} and \textsc{Temam} \cite[Theorem 2.1, p. 243]{EkeTem76CAVP}, there holds 
	\begin{align*}
		\int_\Omega \psi(\xx,\nabla v(\xx))\dd \xx &\leq \liminf_{k\rightarrow \infty}\int_\Omega \psi(\xx,\nabla v_{n_k}(\xx))\dd \xx\\
		& \leq \liminf_{n\rightarrow \infty}\int_\Omega \psi(\xx,\nabla v_{n}(\xx))\dd \xx,
	\end{align*} which shows the lower semicontinuity of $\Psi$ on $\Hi$. \\ Now, we show the \textsc{Lipschitz} continuity of $B$:
	\begin{align*}
		(B(t,u)-B(t,v),w )&= \int_{\Omega}\left( b(t,\xx,u(\xx))- b(t,\xx,v(\xx)) \right) w(\xx) \dd \xx\\
		&\leq \left(\int_{\Omega}\left\vert b(t,\xx,u(\xx))- b(t,\xx,v(\xx)) \right \vert^2 \dd \xx\right)^{1/2}\Vert w\Vert_{\rmL^2(\Omega)}\\
		&\leq \left(\int_{\Omega} C^2_b\vert u(\xx) - v(\xx) \vert^2 \dd \xx\right)^{1/2}\Vert w\Vert_{\rmL^2(\Omega)}  \quad \text{for all }w\in \Hi
	\end{align*} and hence the \textsc{Lipschitz} continuity of $B$. The last assertion follows from Proposition 5.1, p. 21, Proposition 5.7, p. 27, and Proposition 2.1, p. 271, in \textsc{Ekeland} and \textsc{Temam} \cite{EkeTem76CAVP} which completes the proof.
\end{proof}
We are now in the position, to state the main result.
\begin{thm}
	Let the Conditions \ref{2.a}-\ref{2.d} be satisfied and assume that there exists $\tilde{v}\in \rmW^{1,\infty}(\Omega)$ such that $\Psi(\tilde{v})<+\infty$. Then, for every $u_0\in \rmL^2(\Omega) $ and  $v_0\in \overline{\mathrm{dom}(\partial \Psi)}\cap \mathrm{dom}(\Psi)$, there exists a unique weak solution $u\in \rmH^2(0,T,\Hi)$ to \textnormal{(P2)}. Furthermore, the stability estimate \ref{stab} is fulfilled.
\end{thm}
\begin{proof}
	Since the assumptions are all fulfilled, this follows immediately from  Theorem \ref{main}.
\end{proof}
\begin{rem}
	We notice again, that we did not impose any compactness condition on $\Psi$ but rely on the maximal monotonicity of its subdifferential. This implies that we do not identify the nonlinearity in $B$ by compactness but a fixed point argument. Compare with \cite{EmSiWr16ESOT}, where a time discretization method has been employed to show with compactness arguments the existence of weak solutions in an \textsc{Orlicz} space.
\end{rem}

\begin{rem} Instead of the global \textsc{Lipschitz} condition \ref{2.c}, we can impose a local \textsc{Lipschitz} condition on $B$. Consequently, we obtain the local existence result given in Theorem \ref{main}.
\end{rem}

\section{Conclusion}
In this article, we presented a well-posedness result for a doubly nonlinear abstract evolution equation without assuming any linearity of the operators. Our well-posedness proof harnesses a renowned result on contracting semigroups and leverages the \textsc{Banach} fixed-point theorem. Future considerations could delve into more intricate scenarios, such as when the operator $B$ is multi-valued. Further generalizations might pertain to the selection of the underlying space, which could be expanded to a Banach space or even a Gelfand-triplet. However, such extensions would necessitate the employment of advanced techniques and results, akin to those utilized in references like \cite{Akag08DNEB, BaEmMi19EREG} for first-order evolution inclusions.

\bibliographystyle{my_alpha}
\bibliography{bib_aras,alex_pub}

@preamble{
   "\def\cprime{$'$} "}

@article {EmSiWr16ESOO,
    AUTHOR = {Emmrich, Etienne and \v{S}i\v{s}ka, David and Wr\'{o}blewska-Kami\'{n}ska,
              Aneta},
     TITLE = {Equations of second order in time with quasilinear damping:
              existence in {O}rlicz spaces via convergence of a full
              discretisation},
   JOURNAL = {Math. Methods Appl. Sci.},
  FJOURNAL = {Mathematical Methods in the Applied Sciences},
    VOLUME = {39},
      YEAR = {2016},
    NUMBER = {10},
     PAGES = {2449--2460},
      ISSN = {0170-4214},
   MRCLASS = {35L75 (35A01 35L35)},
  MRNUMBER = {3512760},
       DOI = {10.1002/mma.3706},
       URL = {https://doi.org/10.1002/mma.3706},
}

@article {Dion62SPCH,
    AUTHOR = {Dionne, Philippe-A.},
     TITLE = {Sur les probl\`emes de {C}auchy hyperboliques bien pos\'{e}s},
   JOURNAL = {J. Analyse Math.},
  FJOURNAL = {Journal d'Analyse Math\'{e}matique},
    VOLUME = {10},
      YEAR = {1962/63},
     PAGES = {1--90},
      ISSN = {0021-7670},
   MRCLASS = {35.53},
  MRNUMBER = {150475},
MRREVIEWER = {D. Ludwig},
       DOI = {10.1007/BF02790303},
       URL = {https://doi.org/10.1007/BF02790303},
}

@book {Lera53HDE,
    AUTHOR = {Leray, Jean},
     TITLE = {Hyperbolic differential equations},
 PUBLISHER = {The Institute for Advanced Study, Princeton, N. J.},
      YEAR = {1953},
     PAGES = {240},
   MRCLASS = {35.0X},
  MRNUMBER = {0063548},
MRREVIEWER = {H. G. Garnir},
}

@article{Bach22DNEI,
  author       = {Bacho, Aras},
  title        = {Nonsmooth analysis of doubly nonlinear second-order evolution equations with nonconvex energy functionals},
  journal      = {Advances in Nonlinear Analysis},
  year         = {2025},
  volume       = {14},
  number       = {1},
  pages        = {20240064},
  doi          = {10.1515/anona-2024-0064},
  url          = {https://doi.org/10.1515/anona-2024-0064},
}

@article{BACHO2023126,
	author = {Aras Bacho},
	journal = {Journal of Differential Equations},
	pages = {126-169},
	title = {Abstract nonlinear evolution inclusions of second order with applications in visco-elasto-plasticity},
	volume = {363},
	year = {2023}}

@article {Bach23GMYR,
    AUTHOR = {Bacho, Aras},
     TITLE = {A generalization of the {M}oreau-{Y}osida regularization},
   JOURNAL = {J. Math. Anal. Appl.},
  FJOURNAL = {Journal of Mathematical Analysis and Applications},
    VOLUME = {524},
      YEAR = {2023},
    NUMBER = {2},
     PAGES = {Paper No. 127139, 15}
}

@phdthesis {Bach21ONSA,
   author = {Bacho, Aras},
   title = {On the nonsmooth analysis of doubly nonlinear evolution inclusions of first and second order with applications},
   school = {Technische Universit\"{a}t Berlin},
   year = {2021},
   type = {Doctoral Thesis},
   address = {Berlin},
   doi = {10.14279/depositonce-12327},
   url = {http://dx.doi.org/10.14279/depositonce-12327},
}

@article {Akag08DNEB,
    AUTHOR = {Akagi, Goro},
     TITLE = {Doubly nonlinear evolution equations with non-monotone
              perturbations in reflexive {B}anach spaces},
   JOURNAL = {J. Evol. Equ.},
  FJOURNAL = {Journal of Evolution Equations},
    VOLUME = {11},
      YEAR = {2011},
    NUMBER = {1},
     PAGES = {1--41},
}

@article {BaEmMi19EREG,
	AUTHOR = {Bacho, Aras and Emmrich, Etienne and Mielke, Alexander},
	TITLE = {An existence result and evolutionary {$\Gamma$}-convergence
	for perturbed gradient systems},
	JOURNAL = {J. Evol. Equ.},
	FJOURNAL = {Journal of Evolution Equations},
	VOLUME = {19},
	YEAR = {2019},
	NUMBER = {2},
	PAGES = {479--522},
	DOI = {10.1007/s00028-019-00484-x}}

@book {Barb10NDMT,
    AUTHOR = {Barbu, Viorel},
     TITLE = {Nonlinear Differential Equations of Monotone Types in {B}anach Spaces},
 PUBLISHER = {Springer-Verlag},
   ADDRESS = {New York},
      YEAR = {2010},
}

@Book{Barb76NSDE,
  author = 	 {Barbu, V.},
  title = 	 {Nonlinear Semigroups and Differential Equations in 
                  {B}anach spaces},
  publisher = 	 {Noordhoff, Leyden},
  year = 	 {1976},
  OPTannote = 	 {m-accretive},}

@book{Brez73OMMS,
    AUTHOR = {Br{\'e}zis, H.},
     TITLE = {Op{\'e}rateurs Maximaux Monotones et Semi-Groupes de
              Contractions dans les Espaces de {H}ilbert},
 PUBLISHER = {North-Holland Publishing Co.},
   ADDRESS = {Amsterdam},
      YEAR = {1973}}

@book{Brez11FASS,
    AUTHOR = {Br{\'e}zis, Haim},
     TITLE = {Functional Analysis, {S}obolev Spaces and Partial Differential Equations},
 PUBLISHER = {Springer-Verlag},
    ADDRESS = {New York},
      YEAR = {2011}
}

@article {BuMaRa12KVMG,
    AUTHOR = {Bul\'{\i}\v{c}ek, Miroslav and M\'{a}lek, Josef and Rajagopal, K. R.},
     TITLE = {On {K}elvin-{V}oigt model and its generalizations},
   JOURNAL = {Evol. Equ. Control Theory},
  FJOURNAL = {Evolution Equations and Control Theory},
    VOLUME = {1},
      YEAR = {2012},
    NUMBER = {1},
     PAGES = {17--42},
      ISSN = {2163-2472},
       DOI = {10.3934/eect.2012.1.17},
}

@article {BuKaSt1GKVM,
    AUTHOR = {Bul\'{\i}\v{c}ek, Miroslav and Kaplick\'{y}, Petr and Steinhauer, Mark},
     TITLE = {On existence of a classical solution to a generalized
              {K}elvin-{V}oigt model},
   JOURNAL = {Pacific J. Math.},
  FJOURNAL = {Pacific Journal of Mathematics},
    VOLUME = {262},
      YEAR = {2013},
    NUMBER = {1},
     PAGES = {11--33},
      ISSN = {0030-8730},
   MRCLASS = {35Q74 (35A01 35A09 35B40 35B65 74D10 74H20 74H25)},
  MRNUMBER = {3069053},
MRREVIEWER = {Ramon Quintanilla},
       DOI = {10.2140/pjm.2013.262.11},
       URL = {https://doi.org/10.2140/pjm.2013.262.11},
}

@article {CaZh09EUWP,
    AUTHOR = {Cai, Yongyong and Zhou, Shulin},
     TITLE = {Existence and uniqueness of weak solutions for a non-uniformly parabolic equation},
   JOURNAL = {J. Funct. Anal.},
  FJOURNAL = {Journal of Functional Analysis},
    VOLUME = {257},
      YEAR = {2009},
    NUMBER = {10},
     PAGES = {3021--3042},
   
}

@article {Dona74IONP,
    AUTHOR = {Donaldson, Thomas},
     TITLE = {Inhomogeneous {O}rlicz-{S}obolev spaces and nonlinear
              parabolic initial value problems},
   JOURNAL = {J. Differential Equations},
  FJOURNAL = {Journal of Differential Equations},
    VOLUME = {16},
      YEAR = {1974},
     PAGES = {201--256}}

@Book{EkeTem76CAVP,
  AUTHOR = {Ekeland, Ivar and T\'{e}mam, Roger},
     TITLE = {Convex analysis and variational problems},
    SERIES = {Classics in Applied Mathematics},
    VOLUME = {28},
 PUBLISHER = {Society for Industrial and Applied Mathematics (SIAM),
              Philadelphia, PA},
      YEAR = {1999}}

@article {EmmSis11FDSO,
	AUTHOR = {Emmrich, Etienne and {\v{S}}i{\v{s}}ka, David},
	TITLE = {Full discretisation of second-order nonlinear evolution
	equations: strong convergence and applications},
	JOURNAL = {Comput. Methods Appl. Math.},
	VOLUME = {11},
	YEAR = {2011},
	NUMBER = {4},
	PAGES = {441--459},
	DOI = {10.2478/cmam-2011-0025}}

@article {EmmSis13EESO,
	AUTHOR = {Emmrich, Etienne and {\v{S}}i{\v{s}}ka, David},
	TITLE = {Evolution equations of second order with nonconvex potential
	and linear damping: existence via convergence of a full
	discretization},
	JOURNAL = {J. Differential Equations},
	VOLUME = {255},
	YEAR = {2013},
	NUMBER = {10},
	PAGES = {3719--3746},
	DOI = {10.1016/j.jde.2013.07.065}}

@article {EmSiTh15FDNS,
	AUTHOR = {Emmrich, Etienne and {\v{S}}i{\v{s}}ka, David and Thalhammer, Mechthild},
	TITLE = {On a full discretisation for nonlinear second-order evolution
	equations with monotone damping: construction, convergence,
	and error estimates},
	JOURNAL = {Found. Comput. Math.},
	FJOURNAL = {Foundations of Computational Mathematics. The Journal of the
	Society for the Foundations of Computational Mathematics},
	VOLUME = {15},
	YEAR = {2015},
	NUMBER = {6},
	PAGES = {1653--1701},
	DOI = {10.1007/s10208-014-9238-4}}

@article {EmmWro13CQPE,
    AUTHOR = {Emmrich, Etienne and Wr\'{o}blewska-Kami\'{n}ska, Aneta},
     TITLE = {Convergence of a full discretization of quasi-linear parabolic equations in isotropic and anisotropic {O}rlicz spaces},
   JOURNAL = {SIAM J. Numer. Anal.},
  FJOURNAL = {SIAM Journal on Numerical Analysis},
    VOLUME = {51},
      YEAR = {2013},
    NUMBER = {2},
     PAGES = {1163--1184},
}

@article {EmSiWr16ESOT,
	AUTHOR = {Emmrich, Etienne and {\v{S}}i{\v{s}}ka, David and Wr{\'{o}}blewska-Kami{\'{n}}ska,
	Aneta},
	TITLE = {Equations of second order in time with quasilinear damping:
	existence in {O}rlicz spaces via convergence of a full
	discretisation},
	JOURNAL = {Math. Methods Appl. Sci.},
	VOLUME = {39},
	YEAR = {2016},
	NUMBER = {10},
	PAGES = {2449--2460},
	DOI = {10.1002/mma.3706}}

@article {EmmTha11DNEE,
	AUTHOR = {Emmrich, Etienne and Thalhammer, Mechthild},
	TITLE = {Doubly nonlinear evolution equations of second order:
	existence and fully discrete approximation},
	JOURNAL = {J. Differential Equations},
	VOLUME = {251},
	YEAR = {2011},
	NUMBER = {1},
	PAGES = {82--118},
	DOI = {10.1016/j.jde.2011.02.014}}

@article {EmmTha10CTDD,
	AUTHOR = {Emmrich, Etienne and Thalhammer, Mechthild},
	TITLE = {Convergence of a time discretisation for doubly nonlinear
	evolution equations of second order},
	JOURNAL = {Found. Comput. Math.},
	VOLUME = {10},
	YEAR = {2010},
	NUMBER = {2},
	PAGES = {171--190},
	DOI = {10.1007/s10208-010-9061-5}}

@article {FriNec88SNWE,
    AUTHOR = {Friedman, Avner and Ne\v{c}as, Jind\v{r}ich},
     TITLE = {Systems of nonlinear wave equations with nonlinear viscosity},
   JOURNAL = {Pacific J. Math.},
  FJOURNAL = {Pacific Journal of Mathematics},
    VOLUME = {135},
      YEAR = {1988},
    NUMBER = {1},
     PAGES = {29--55},
      ISSN = {0030-8730},
   MRCLASS = {35L70 (35K55 73F99)},
  MRNUMBER = {965683},
MRREVIEWER = {Reza Malek-Madani},
       URL = {http://projecteuclid.org/euclid.pjm/1102688343},
}

@incollection {GwSw11PEAO,
    AUTHOR = {Gwiazda, Piotr and \'{S}wierczewska-Gwiazda, Agnieszka},
     TITLE = {Parabolic equations in anisotropic Orlicz spaces with
              general $\mathscr{N}$-functions},
 BOOKTITLE = {Parabolic problems},
    SERIES = {Progr. Nonlinear Differential Equations Appl.},
    VOLUME = {80},
     PAGES = {301--311},
 PUBLISHER = {Birkh\"{a}user/Springer Basel AG, Basel},
      YEAR = {2011}}

@book {Lion69QMRP,
    AUTHOR = {Lions, J.-L.},
     TITLE = {Quelques M\'{e}thodes de R\'{e}solution des Probl\`emes aux Limites Non Lin\'{e}aires},
 PUBLISHER = {Dunod; Gauthier-Villars},
   ADDRESS = {Paris},
      YEAR = {1969}}

@article {LioStr65SNLE,
    AUTHOR = {Lions, Jacques-Louis and Strauss, W. A.},
     TITLE = {Some non-linear evolution equations},
   JOURNAL = {Bull. Soc. Math. France},
  FJOURNAL = {Bulletin de la Soci\'{e}t\'{e} Math\'{e}matique de France},
    VOLUME = {93},
      YEAR = {1965},
     PAGES = {43--96},
      ISSN = {0037-9484},
   MRCLASS = {34.95},
  MRNUMBER = {199519},
MRREVIEWER = {T. Kato},
       URL = {http://www.numdam.org/item?id=BSMF_1965__93__43_0},
}

@book{LioMag72NHBV,
    AUTHOR = {Lions, J.-L. and Magenes, E.},
     TITLE = {Probl\`emes aux Limites Non Homog\`enes et Applications. {V}ol. 1},
    SERIES = {Travaux et Recherches Math\'{e}matiques, No. 17},
 PUBLISHER = {Dunod, Paris},
      YEAR = {1968}}

@article {Puhs15OEVF,
    AUTHOR = {Puhst, Dimitri},
     TITLE = {On the evolutionary fractional {$p$}-{L}aplacian},
   JOURNAL = {Appl. Math. Res. Express. AMRX},
  FJOURNAL = {Applied Mathematics Research Express. AMRX},
      YEAR = {2015},
    NUMBER = {2},
     PAGES = {253--273},
      ISSN = {1687-1200},
   MRCLASS = {35R11},
  MRNUMBER = {3394266},
       DOI = {10.1093/amrx/abv003},
       URL = {https://doi.org/10.1093/amrx/abv003},
}

@article {RosTho17CDIP,
    AUTHOR = {Rossi, Riccarda and Thomas, Marita},
     TITLE = {Coupling rate-independent and rate-dependent processes:
              existence results},
   JOURNAL = {SIAM J. Math. Anal.},
  FJOURNAL = {SIAM Journal on Mathematical Analysis},
    VOLUME = {49},
      YEAR = {2017},
    NUMBER = {2},
     PAGES = {1419--1494},
       DOI = {10.1137/15M1051567},
}

@Book{Roub13NPDE,
  AUTHOR = {Roub\'{\i}\v{c}ek, Tom\'{a}\v{s}},
     TITLE = {Nonlinear Partial Differential Equations with Applications},
    SERIES = {International Series of Numerical Mathematics},
    VOLUME = {153},
   EDITION = {2nd},
 PUBLISHER = {Birkh\"{a}user, Basel},
      YEAR = {2013}}

@article {Ruf17CFDS,
    AUTHOR = {Ruf, A. M.},
     TITLE = {Convergence of a full discretization for a second-order
              nonlinear elastodynamic equation in isotropic and anisotropic
              {O}rlicz spaces},
   JOURNAL = {Z. Angew. Math. Phys.},
  FJOURNAL = {Zeitschrift f\"{u}r Angewandte Mathematik und Physik. ZAMP.
              Journal of Applied Mathematics and Physics. Journal de
              Math\'{e}matiques et de Physique Appliqu\'{e}es},
    VOLUME = {68},
      YEAR = {2017},
    NUMBER = {5},
     PAGES = {Paper No. 118, 24},
      ISSN = {0044-2275},
   MRCLASS = {35L75 (35L35 35L82 47H05 47J35 65M12 65M60)},
  MRNUMBER = {3705806},
       DOI = {10.1007/s00033-017-0863-z},
       URL = {https://doi.org/10.1007/s00033-017-0863-z},
}

@book{Wlok82PD,
author =  {Wloka, Joseph},
title =   {Partielle {D}ifferentialgleichungen},
year = "1982",
address = {Leipzig},
publisher = "B. G.~Teubner"}

@Book{Zeid90NFA2b,
    AUTHOR = {Zeidler, Eberhard},
     TITLE = {Nonlinear Functional Analysis and its Applications {II}/{B}},
 PUBLISHER = {Springer-Verlag},
   ADDRESS = {New York},
      YEAR = {1990}}

@Book{Zeid85NFA3,
    AUTHOR = {Zeidler, Eberhard},
     TITLE = {Nonlinear functional analysis and its applications {I}{I}{I}},
      NOTE = {Variational methods and optimization},
 PUBLISHER = {Springer-Verlag},
   ADDRESS = {New York},
      YEAR = {1985},
  OPTPAGES = {xxii+662},}

%
%

%

\end{document}